\documentstyle{amsppt}
\NoRunningHeads
\NoBlackBoxes
\topmatter
\title Great sphere foliations and manifolds with curvature bounded above
\endtitle
\author Vladimir ROVENSKII and Victor TOPONOGOV\endauthor
\address Chair of Geometry, Mathematical Dept.,
Pedagogical State University,
\newline Lebedeva st. 89, Krasnoyarsk-49, 660049, Russia\endaddress
\email roven\,\@\,edk.krasnoyarsk.su\endemail
\address Mathematical Institute, 630090 Novosibirsk--90, Russia
\endaddress
\subjclass 53C12, 53C20\endsubjclass
\thanks This paper is a final form and no version of it will be submitted for
publication elsewhere.\endthanks
\abstract The survey is devoted to Toponogov's conjecture, that {\it if a
complete simply connected Riemannian manifold with sectional curvature $\le 4$
and injectivity radius $\ge\frac\pi 2$ has extremal diameter $\frac\pi 2$,
then it is isometric to CROSS}. In Section 1 the relations of problem with
geodesic foliations of a round sphere are considered, but the proof of
conjecture on this way is not complete. In Section 2 the proof based on recent
results and methods for topology and volume of Blaschke manifolds is given.
\endabstract
\endtopmatter
\document
\head 1. Great circle foliations\endhead

 The interest to fibrations of the $n-$sphere by great $\nu-$spheres is
propelled by the Blaschke problem and by extremal theorems in Riemannian
geometry.

On a round sphere a foliation by geodesics is the same thing as a great circle
fibration. The most simple of them -- {\it Hopf fibration\/} with fibers
$\{S^1\}$, can be given as a collection of intersections of $S^{2n-1}$ with
all holomorphic $2-$planes $\{\sigma=x\land Jx\}$, where $J$ is a {\sl complex
structure} -- a linear operator in $\Bbb R^{2n}$, given for some orthonormal
basis $\{e_i\}$ by the rule
$$Je_{2i-1}=e_{2i},\ Je_{2i}=-e_{2i-1},\quad(1\le i\le n).\eqno(1)$$
 Let $\Cal F_0(S^{2n-1})$ denote a {\sl space of all Hopf fibrations} of
sphere $S^{2n-1}$, for $n=2$ see Chapter 1. Each fiber spans corresponding
oriented two-plane through the origin in $\Bbb R^{2n}$ and hence determines a
point in the Grassmann manifold $G(2,2n)$.

\definition{2 Definition ([GW] for n=2)} The {\it skew-Hopf fibration\/}
is given by intersections of $S^{2n-1}$ with all holomorphic
$2-$planes $\{\sigma=x\land Jx\}$, where $J$ is an {\sl almost complex
structure} -- linear operator in $\Bbb R^{2n}$, given for some affine basis
$\{e_i\}$ by the rule (1).\enddefinition

 By other words, the skew-Hopf fibration is obtained from the Hopf fibration
by applying a nondegenerate linear transformation of $\Bbb R^{2n}$ and then
projection the images of fibers back to the $S^{2n-1}.$

 Let $\Cal F_1(S^{2n-1})$ denote a {\sl space of all skew-Hopf fibrations} of
sphere $S^{2n-1}$.
Let $\Cal F(S^{2n-1})$ denote a {\sl space of all oriented great circle
fibrations} of sphere $S^{2n-1}$.

\vskip 5pt
For simplicity we consider below geodesic foliations on three-sphere
(for $n>2$ see [Rov 2]).
 The space $\Cal F_0(S^3)$ is $2-$dimensional and homeomorphic to a pair of
disjoint two-spheres. The $8-$dimensional space $\Cal F_1(S^3)$ is a disjoint
union of two copies of $S^2\times\Bbb R^6$. Both these spaces are
homogeneous [GW]:
$$\Cal F_0(S^3)=
O(4)/U(2),\ \Cal F_1(S^3)=GL(4,\Bbb R)/GL(2,\Bbb C).\eqno(3)$$
The space $\Cal F(S^3)$ is infinite dimensional. Let $V$ be a Hopf unit vector
field and $D^2$ -- a small ball transverse to $V$ at point $p$. As was shown
in [GW], a small $C^1-$perturbations of $V$ on $D^2$, which are identity on
the neighbourhood of boundary $\partial D^2$, lead to different great circle
foliations of $S^3$.

 For example, $G(2,4)\equiv S^2\times S^2$ [GW], and a Hopf fibration
$h\in\Cal F_0(S^3)$ can be recognized by the fact that its orbit space $M_h$
appears inside the Grassmannian as $\{point\}\times S^2$ or
$S^2\times\{point\}$.

For two $2-$planes $P$ and $Q$ in $\Bbb R^4$ the smallest and the largest
angles $\alpha_{min}$ and $\alpha_{max}$ from interval $[0,\frac\pi 2]$, that
any line in $P$ makes with the plane $Q$, are called {\it the principal angles}
between certain planes. The same angles result upon interchanging the roles of
$P$ and $Q$. The relative position of $P$ and $Q$ in $\Bbb R^4$ is completely
determined by these principal angles in sence of rigid motion of $\Bbb R^4$.
One can always choose an orthonormal basis $e_1,e_2,e_3,e_4$ for $\Bbb R^4$ so
that $e_1,e_2$ is an orthonormal basis for $P$ and
$\cos\alpha_{min}e_1+\sin\alpha_{min}e_3,
\ \cos\alpha_{max}e_2+\sin\alpha_{max}e_4$
is an orthonormal basis for $Q$ [GW]. Any two equidistant great circles from
$S^3$ (in particular two leaves of Hopf fibration) determine two $2-$planes
from $\Bbb R^4$ with equal principal angles.

 Not all submanifolds of $G(2,4)$ can appear in the role of $M_f$ for
$f\in\Cal F(S^3)$.

\proclaim{4 Theorem [GW]} A submanifold in $G(2,4)\equiv S^2\times S^2$
corresponds to a fibration $f\in\Cal F(S^3)$ if and only if it is a graph of a
distance decreasing map $\tilde f$ from either $S^2$ factor to the other.
Fibration $f$ is differentiable if and only if the corresponding map $\tilde
f$ is differentiable with $|df|\le 1$.\endproclaim

 Hence, the space $\Cal F(S^3)$ deformation retracts to the space
$\Cal F_0(S^3)$.

 The catalogue of great circle fibrations of the three-sphere (in Theorem
4) is one of first nontrivial examples in which one has a clear overview of
all possible geodesic foliations of a fixed Riemannian manifold.
There are plenty of nondifferentiable great circle fibrations of $S^3$. There
also exist discontinuous fillings of $S^3$ by great circles: one can fill the
closed solid torus $x_1^2+x_2^2\ge x_3^2+x_4^2$ on $S^3$ as for Hopf original
fibration and then fill the remaining open solid torus as for Hopf fibration
with reversed screw sence [GW].

\proclaim{5 Corollary [GW]} Each fibration $f\in\Cal F(S^3)$ contains a
pair of orthogonal fibers.\endproclaim

\proclaim{6 Theorem [Gag]} Each skew-Hopf fibration $f\in\Cal F_1(S^3)$
corresponds to a distance decreasing map $\tilde f:S^2\to S^2$ with convex
image in a semi-sphere, which can be decomposed as: \ a) orthogonal projection
of $ S^2$ to a plane through the center of the sphere; followed by \ b) a
distance decreasing linear map from one 2-plane to another; and finally
\ c) inverse projection onto $ S^2$.\endproclaim

 Thus the space $\Cal F_1(S^3)$ is stratified with respect to values of rank
$r$ of linear map from point b):
each fibration with $r=2$ has exactly one pair of antipodes, in degenerate
case $r=1$ each fibration has $1-$parameter family (Hopf torus) of pairs of
antipodes, case $r=0$ is the same as Hopf fibrations.

 By Corollary 5 for each $f\in\Cal F(S^3)$ there exists orthonormal {\it
special basis} $\{e_i\}$ in $\Bbb R^4$ such, that the sphere $S^3$ intersects
with $2-$planes $e_1\land e_2$ and $e_3\land e_4$ by pair of orthogonal
fibers, (which is unique when $r=2$).
It is not difficult to see, that in special basis for each vector
$\lambda=(1,0,\lambda_3,\lambda_4)$ there exists a unique vector $
h=(0,1,h_3,h_4)$ such, that the intersection of the plane $\lambda\land h$
with sphere $ S^3$ is a fiber of $f$.  Hence, we may correspond to foliation
$f\in\Cal F(S^3)$ a diffeomorphism of $2-$plane
$$\varphi=(\varphi_1,\varphi_2):\Bbb R^2\to\Bbb R^2,
\ \ h_3=\varphi_1(\lambda_3,\lambda_4),\ h_4=\varphi_2(\lambda_3,\lambda_4).
\eqno(7)$$
 If $f\in\Cal F_0(S^3)$, then the corresponding $\varphi:\Bbb R^2\to\Bbb R^2$
is a linear orthogonal map.

\proclaim{8 Lemma [Rov 1]} The space $\Cal F_1(S^3)$ is characterized in
$\Cal F(S^3)$ by the property, that $\varphi:\Bbb R^2\to\Bbb R^2$ is a linear
operator and without real eigenvalues.\endproclaim

\demo{Proof} Let $f\in\Cal F_1(S^3)$, and an almost complex structure $J$ in
$\Bbb R^4$ for special basis is given by block matrix:
$$J=\left(\matrix A&0\\
0&B\endmatrix\right),
\ A=\left(\matrix a_{11}&a_{12}\\
a_{21}&a_{22}\endmatrix\right),
\ B=\left(\matrix a_{33}&a_{34}\\
a_{43}&a_{44}\endmatrix\right),\ A^2=B^2=-E,\eqno(9)$$
moreover, the matrices $A,\,B$ have no real eigenvalues. It is easy to see,
that $\varphi$ is a linear operator on $Bbb R^2$ given by the following matrix
without real eigenvalues
$$F=\left(\matrix a&b\\
c&d\endmatrix\right),
\ a=\frac{a_{33}-a_{11}}{a_{21}},\ b=\frac{a_{34}}{a_{21}},
\ c=\frac{a_{43}}{a_{21}},\ d=\frac{a_{44}-a_{11}}{a_{21}}.\eqno(10)$$

 Conversely, let $\varphi$ be a linear transformation of $\Bbb R^2$ with
matrix
$F=\left(\smallmatrix a&b\\ c&d\endsmallmatrix\right)$
without real eigenvalues, i.e. the discriminant of characteristic quadratic
equation is negative $ D=(a-d)^2+4bc<0.$
Then the matrix of form (9) with coefficients
$$a_{11}=-a_{22}=-\frac{a+d}{\sqrt{-D}},
\ a_{33}=-a_{44}=\frac{a-d}{\sqrt{-D}},$$
$$a_{12}=\frac{2(bc-ad)}{\sqrt{-D}},\ a_{21}=\frac 2{\sqrt{-D}},
\ a_{34}=\frac{2b}{\sqrt{-D}},\ a_{43}=\frac{2c}{\sqrt{-D}}\eqno(11)$$
defines an almost complex structure in $\Bbb R^4$. It is easy to see, that
the for induced foliation $\tilde f\in\Cal F_1(S^3)$ the corresponding
operator $\tilde\varphi:\Bbb R^2\to\Bbb R^2$ has the following matrix
$F=\left(\smallmatrix a&b\\ c&d\endsmallmatrix\right)$.\qed\enddemo

 We study the skew-Hopf fibrations in relations with interesting space
$\Cal F_R(S^{2n-1})$ of analytic geodesic foliations of sphere $S^{2n-1}$,
see [Top 1-3] for more general case.

\definition{12 Definition [Top 1]} The subspace
$\Cal F_R(S^{2n-1})\subset\Cal F(S^{2n-1})$
consists of fibrations, for which there exists a tensor (multi-linear
function)
$R:\Bbb R^4\times\Bbb R^4\times\Bbb R^4\times\Bbb R^4\to\Bbb R^1$
with properties:

$R_1)$ curvature symmetries
$$R(x,y,z,w)=R(z,w,x,y)=-R(x,y,w,z),$$
$$R(x,y,z,w)+R(z,x,y,w)+R(y,z,x,w)=0$$

$R_2)$ for almost each unit vector $x\in\Bbb R^{2n}$ there exists unique
$2-$plane $\sigma\owns x,$ with condition
$$R(x,y,x,y)=1,\ (y\in\sigma,\ y\ \bot\ x,\ |y|=1),$$

$R_3)$ if a $2-$plain $\sigma$ contains a fiber of $f$, then
$$R(x,y,x,y)=1,\ (x,y\in\sigma,\ y\ \bot\ x,\ |x|=|y|=1),$$
i.e. the {\it sectional curvature} of $R$ is extremal on $\sigma$.
\enddefinition

 Note, that the curvature tensor on tangent spaces of $CP^n$ with complex
structure $J$, standard metric and corresponding Hopf fibration have the
properties $R_1)-R_3)$.

\proclaim{13 Theorem [Rov 1]} $\Cal F_R(S^3)=\Cal F_1(S^3)$.\endproclaim

\remark{14 Remark} The analogous fact is true for $\Cal F_R(S^{2n-1})$ with
$n>2$, [Rov 2]. From theorem 13 it follows the negative answer for
assumption by [Lemma 6, Top 1] of equality $\Cal F_R(S^3)=\Cal F_0(S^3)$, also
see 2.\endremark

\demo{Proof} For each $f\in\Cal F(S^3)$ there exists special orthonormal basis
such, that $2-$pla\-ins $e_1\land e_2$ and $e_3\land e_4$ intersect with
sphere $S^3$ by fibers. Let $\varphi:\Bbb R^2\to\Bbb R^2$ be a diffeomorphism,
which by Lemma 8 corresponds to fibration $f\in\Cal F(S^3)$.

1. We shall state firstly the inclusion $\Cal F_R(S^3)\subset\Cal F_1(S^3)$.

Let $f\in\Cal F_R(S^3)$. For a special basis in view of $R_1)-R_3)$ we have
$$R_{121i}=\delta_{i2},\ R_{212i}=\delta_{i1},
\ R_{343i}=\delta_{i4},\ R_{434i}=\delta_{i3},$$
where $R_{ijkp}=R(e_i,e_j,e_k,e_p)$ are components of tensor $R$.
Let us calculate the {\it sectional curvature} of $R$:
$$K(\lambda,h)=\frac{R(\lambda,h,\lambda,h)}{\lambda^2h^2-(\lambda,h)^2}=
\frac Q{\lambda^2h^2-(\lambda,h)^2}+1.$$
 The quadratic form $Q$ from variables $\lambda_3,\lambda_4,h_3,h_4$ is given
by formula
$$Q=\sum_{i,k}(R_{2i2k}-\delta_{ik})\lambda_i\lambda_k+\sum_{j,p}(R_{1j1p}-
\delta_{jp})h_jh_p+2\sum_{i,p}(R_{12ij}+R_{2ij1})\lambda_ih_j.\eqno(15)$$
The $2-$parameter family of vectors $(\lambda_3,\lambda_4,h_3,h_4)$, which
correspond to fibers of $f$, in view of $R_2)-R_3)$ lies in kernel of form
$Q$. In view of $R_3)$ the subspace $\ker Q$ is $2-$dimensional, and hence the
functions $h_3(\lambda_3,\lambda_4),\ h_4(\lambda_3,\lambda_4)$ of variables
$\lambda_3,\ \lambda_4$ are linear. Since the fibers do not intersect in
$S^3$, then linear operator $\varphi$  has no real eigenvalues.
By Lemma 8 $f\in\Cal F_1(S^3)$.

2. We shall now state the inverse inclusion
$\Cal F_1(S^3)\subset\Cal F_R(S^3)$.

Let $f\in\Cal F_1(S^3)$. Denote by
$F=\left(\smallmatrix a&b\\ c&d\endsmallmatrix\right)$
the matrix of linear operator $\varphi$, which corresponds to $f$ in special
ortho-basis, see Lemma 8. In view multi-linearity and symmetries of tensor
$R$ it is sufficient to define only some of components $\{R_{ijkp}\}$.
It is easy to see, that properties $R_1)-R_3)$ for tensor $R$ and $f$ are
true, when we assume for arbitrary $\gamma<0,\ \beta<0$
$$R_{121i}=\delta_{i2},\ R_{2i2i}=\delta_{i1},\ R_{343i}=\delta_{i4},
\ R_{434i}=\delta_{i3},\ R_{1313}-1=\gamma,$$
$$R_{1414}-1=\beta,\ R_{2424}-1=b^2\gamma+d^2\beta,\ R_{1323}=a\gamma,
\ R_{2324}=ab\gamma+cd\beta,$$
$$R_{1234}=\frac 13(-c\beta+b\gamma),\ R_{2323}-1=a^2\gamma+c^2\beta,
\ R_{1424}=d\beta,\ R_{1314}=0,$$
$$R_{2431}=-\frac 13(2b\gamma+c\beta),\ R_{2341}=-\frac 13(2c\beta+b\gamma).
\eqno(16)$$
Thus $f\in\Cal F_R(S^3)$.\qed\enddemo

\head 2. Extremal theorem for manifolds with curvature bounded above\endhead
\rightheadtext{Manifolds with curvature bounded above}

 A complete even-dimensional simply connected Riemannian manifold $M$ with
sectional curvature $0<K_M\le 4$ has injectivity radius
$inj\,(M)\ge\frac\pi 2$ and, hence, the diameter $diam\,(M)\ge\frac\pi 2$.
 For complete odd-dimensional simply connected $M$ the same inequalities for
$inj\,(M)$ and $diam\,(M)$ are true under more strong curvature restriction
$1\le K_M\le 4$. The optimal value for the curvature pinching constant
$\delta$ in the last case is unknown [AM], but for $\delta<\frac 19$
the proposition is wrong in view of Berge's example [Ber 2].

\remark{17 Remark} M.\,Berge considered a family of Riemannian metrics
$g_s, (0<s\le 1)$ on the odd-dimensional spheres $S^{2n+1}$, which are
defined by shrinking the standard metric in the direction of the Hopf circles
in such a way that their lengths with respect to $g_s$ become $2\pi s$.
The range of the sectional curvature of such a metric $g_s$ is in the interval
$[s^2,4-3s^2]$. Clearly, $\pi s<\pi/\sqrt{4-3s^2}$ for $s^2<\frac 13$. This
means that for any $\delta\in(0,\frac 19)$ there exists a M.\,Berger metric
$g_s$ whose sectional curvature $K$ is $\delta-$pinched and whose injectivity
radius is strictly less than $\pi/\sqrt{\max K}$.
>From considering the curvature along horizontal geodesics, we obtain that for
any $0<s<1$ the conjugate radius of $g_s$ is strictly greater than
$\pi/\sqrt{\max K}$.\endremark

 Note, that $inj\,(M)$ is always not more then $diam\,(M)$. The manifolds with
$inj\,(M)=diam\,(M)$ are exactly Blaschke manifolds [Bes].

 In situations when the extremal value for curvature, diameter or volume of
manifold is considered (under others given conditions), one often obtains,
that manifold is isometric to a model space from a finite list. The bright
example of such {\it extremal theorems} is the following

\proclaim {18 Theorem (minimal diameter) [Ber 1]} Let $M$ be a complete,
connected, simply connected Riemannian manifold with sectional curvature
$1\le K_M\le 4$ and diameter $\frac\pi 2$. Then $M$ is isometric to CROSS:
round sphere of curvature $4$ or projective spaces $CP^n,\,HP^n,\,CaP^2$ with
its canonical metric.\endproclaim

 M.\,Berger used direct geometric arguments to see the curvature tensor in
enough detail to prove that such a manifold must be locally symmetric and
hence (since simply connected) symmetric space. Appealing to the
classification of these finished the proof. J.\,Cheeger and D.\,Ebin [CE] also
give a geometricaly more direct proof for this result. H.\,Gluck and
co-authors [GWZ] reproved constructively this theorem, used Berger's geometric
arguments to show that the exponential map from the tangent cut locus to the
cut locus is a fibration of a round sphere by parallel great spheres, and
hence a Hopf fibration. Then they see how this fibration encodes the curvature
tensor and use this to display an isometry between $M$ and a round sphere or
projective space. Latter Berger's theorem was generalized in some directions:

1) stability results:
{\sl there exists a constant $\delta_n<\frac 14$ such that any
$n-$dimensional, complete, simply connected Riemannian manifold $M^n$ with
$\delta_n\le K_M\le 1$ is homeomorphic to CROSS} [Ber 3],
{\sl for compact odd-dimensional manifold such $\delta$ is universal and less
than} $\frac 14(1+10^{-6})^{-2}$ [AM].

2) upper curvature bound is replaced by corresponding lower bound for the
diameter or radius:
{\sl a complete, simply connected Riemannian manifold $M^n$ with $K_M\ge 1$
and radius $rad\,(M)\ge\frac\pi 2$ is either homeomorphic to the sphere or
the universal covering $\tilde M^n$ is isometric to CROSS} [Wil].

Analogous result with diameter is obtained by [GG3], where the only case of
$\Bbb CaP^2$ is unknown.
Recall that the {\it radius} $rad\,(M)$ of a compact connected Riemannian
manifold is defined as the infimum of the function
$p\to rad\,_p(M)=:\max_{q\in M}dist\,(p,q)$. Clearly,
$inj\,(M)\le rad\,(M)\le diam\,(M)$.
One of key point in these results is the studying of Riemannian foliations on
the round sphere. For 1- and 3- dimensional leaves they are always Hopf
fibrations [GG], a partial classification of Riemannian foliations on
$S^{15}$ with 7-dimensional leaves is obtained by [Wil], [Lu].

\vskip 10pt
 Thus it is natural to investigate $V^m(-\infty,4)$ -- a complete simply
connected Riemannian manifold, whose sectional curvature $K_V\le 4$ and
injectivity radius $inj\,(V)\ge\frac\pi 2$. Since this manifold has lower
estimate for diameter, the case of extremal value $\frac\pi 2$ of diameter
is especially interesting.
Note that inequality $inj\,(V)\ge\frac\pi 2$ is equivalent to condition, that

{\sl the perimeter of every nondegenerate geodesic biangle in $V^m(-\infty,4)$
is not less than $\pi$}.

\proclaim {19 Theorem [Top 2,3]} The manifold $M=V^{2n+1}(-\infty,4)$ with
diameter $\frac\pi 2$ is isometric to sphere of curvature $4$.
The manifold $M=V^{2n}(-\infty,4)$ with extremal diameter $\frac\pi 2$ is
isometric to sphere of curvature $4$, either geodesics are grouped into
families (as for projective spaces):

$F_1)$ for every point $p\in M$ and any vector $\lambda\in T_pM$ there exists
$a-$dimensional ($a=2,\,4,\,8$ and if $a=8$, then $\dim M=16$) subspace
$d(\lambda)\subset T_pM$ such, that all geodesics
$\gamma\subset M,\ (\gamma(0)=p,\ \gamma'(0)\in d(\lambda))$ form a totally
geodesic submanifold $F(p,\lambda)$, which is isometric to round sphere
$S^a(4)$;

$F_2)$ for all nonzero vectors $\lambda_1,\,\lambda_2\in T_pM$ the
submanifolds $F(p,\lambda_1),\,F(p,\lambda_2)$ coincide, either their
intersection consists only of one point $p$.\endproclaim

The key point of Theorem 19 is the following result:
\proclaim {20 Theorem [Top 2]} If a Riemannian manifold $V^n(-\infty,4)$ has a
closed geodesic $\gamma$ with length $\pi$ and index $a-1$, then there is
$a-$dimensional totally geodesic submanifold containing $\gamma$, which is
isometric to $a-$dimensional sphere of curvature~$4$.\endproclaim

 We outline the idea of the proof of Theorem 20.
\proclaim {21 Lemma [Top 2]} Under the conditions and with the notations of
Theorem 20, any two points $P$ and $Q$ of $\gamma$ whose mutual distance on
$\gamma$ is $\frac\pi 2$ are conjugate of multiplicity $a-1$.\endproclaim

 The proof of Lemma 21 uses the condition $inj\,(V)\ge\frac\pi 2$ and
well-known Lemma of the calculus of variations.

 From Lemma 21 we obtain by an easy induction
\proclaim {22 Lemma [Top 2]} Under the conditions and with the notations
of Theorem 20, there exists an arc $\sigma$ of length $>\frac\pi 2$ on the
geodesic $\gamma$ and a $(a-1)-$parameter family of parallel vector fields
$\nu$ along $\sigma$ such that the Riemann curvature for the two-dimensional
directions in $\nu$ along $\gamma$ is equal to $4$.\endproclaim

 Using Lemma 22 for all the fields $\nu$ we can construct a sequence of
triangles $\Delta_n(\nu)$ whose perimeter is strictly less then $\pi$ and
which converges to $\gamma$. It follows from condition $inj\,(V)\ge\frac\pi 2$
that in every triangle $\Delta_n(\nu)$ we can span a cone $K_n(\nu)$ obtained
as the set of the shortest lines between the vertices of $\Delta_n(\nu)$ and
the opposite edges. For these cones the following Lemma holds.

\proclaim {23 Lemma [Top 2]} The Gauss curvature in points of $K_n(\nu)$ does
not exceed $4$.\endproclaim
 Lemma 23 follows from Synge's Theorem.

\proclaim {24 Lemma [Top 2]} The area of the cone $K_n(\nu)$ is not greater
than the area of the triangle $\Delta_n^L(\nu)$ on the sphere of curvature $4$
whose sides have the length of the corresponding sides of $\Delta_n(\nu)$.
\endproclaim
 Lemma 24 follows from a Theorem of A.D.Aleksandrov, see [Top 2]. From Lemmas
23 and 24 we obtain an upper bound of the integral curvature of the cones
$K_n(\nu)$. On the other hand, the Gauss-Bonnet Theorem for the integral
curvature of $K_n(\nu)$ can be used to obtain a lower bound for the angle sum
of $K_n(\nu)$.
 A comparison of these bounds shows that the Gauss curvature of $K_n(\nu)$ is
everywhere almost equal to $4$ and the area of $K_n(\nu)$ is almost equal to
$\pi$. Passing to the limit for $n\to\infty$, we see that there exists a
$(a-1)-$parameter family of surfaces $\{F\}$ that are isometric to the
$2-$dimensional hemisphere of curvature $4$ and whose boundary is $\gamma$.

 It is now easily shown that the union of all surfaces of that family is a
$a-$dimensi\-onal surface $F^a$ which is isometric to the $a-$dimensional
sphere. From the previous results and $inj\,(V)\ge\frac\pi 2$ it follows
easily that $F^a$ is a totally geodesic surface in $M$.

 The proof of Theorem 19 follows similar reasoning, analog to the preceding
reduction, only we need the certain topological results and in particular a
Theorem of W.Brouder.

\vskip 3mm
V.\,Toponogov [Top 1-3] {\bf conjectured}, that {\it a manifold
$V^{2n}(-\infty,4)$ with extremal diameter $diam\,(V)=\frac\pi 2$ is isometric
to CROSS}.

Note, that the tangent $a-$planes to submanifolds $\{F(p,*)\}$ (in Theorem
19) induce $(a-1)-$dimensional great sphere foliation of the round sphere
$S_p$ in the tangent space $T_pM$. In case $K_M>0$ for almost each point
$p\in M$ such foliation $f_p$ is related with function of sectional curvature
in $T_pM$ by the following way (compare with 12) [Top 1]:

$\bar R_2)$ {\sl for almost each vector $x\in T_pM$ there exists unique
$a-$dimensional subspace $V\owns x$ with condition}
$$K(x,y)=4,\ (y\in V),$$

$\bar R_3)$ {\sl if $a-$dimensional subspace $V$ contains a fiber of $f_p$,
then}
$$K(x,y)=4,\ (x,\,y\in V).$$

 The natural strategy to prove, that {\sl Riemannian manifold $M$ with
properties $F_1),\,F_2)$ is isometric to CROSS,} (and, hence, to prove
Toponogov's conjecture) is to deduce firstly, that induced great sphere
foliations in tangent spheres $\{S_p\},\ (p\in M)$ are Hopf fibrations.
 The last claim was conjectured in [Lemma 6, Top 1] for foliation with
properties $\bar R_2),\,\bar R_3)$ and curvature symmetries $R_1)$.
>From the consideration at one point of $M$ (see results of 1) one can only
obtain, that such foliations on tangent spheres $\{S_p\},\ (p\in M)$, are
skew-Hopf fibrations (the proof of [Lemma 6, Top 1] is wrong).
By other words, manifold in Theorem 19, when $a=2$, admits an almost complex
structure $J:TM\to TM$ with identity
$$(\nabla_xJ)x=0,\ (x\in TM),$$
and with constant holomorphic curvature (totally geodesic submanifolds
$F(p,\lambda)$ are $J-$invariant). Note, that Hermitian manifold $(M,\,J)$
with above identity is called {\it nearly Kahlerian}. We don't know is it
possible to prove locally, that $J$ is Hermitian.

\vskip 10pt
 Latter the Toponogov's conjecture was proved on another way [RovT]:
by using of global integral geometrical methods by M.\,Berge and J.\,Kazdan
and recent topological results for manifolds with closed geodesics.

 Manifolds $M$ with properties $F_1),\,F_2)$ are the particular case of
Blaschke manifolds. We shall give a short survey.

\vskip 10pt
{\bf 25. Manifolds with closed geodesics}.
A compact Riemannian manifold $M$ is called a $C_\pi-${\it manifold}, if all
its geodesics are closed of equal length $\pi$. This class includes
(A.Allamigeon and F.Warner) {\it Blaschke manifolds}, for which all cut loci
$Cut(p)\subset T_pM,\ (p\in M)$ are round spheres of constant radius and
dimension. The examples are CROSS: a sphere or projective spaces over a
classical fields.

 If $M$ is a simply connected $C_\pi-$manifold, then it is homotopically
equivalent to CROSS (R.Bott and H.Samelson).

For Blaschke manifold the exponential map $\exp_p:T_pM\to C(p)$, restricted on
sphere $S_p$ with radius $d(p,C(p))$, defines a {\it great sphere foliation},
for CROSS this foliation is Hopf fibration. Since every great sphere foliation
of $S^N$ is homeomorphic to Hopf fibration [Sat], (partial cases in works by
H.Gluck, F.Warner, C.Yang), then a simply connected Blaschke manifold is
homeomorphic to its model CROSS.

 The well-known {\bf Blaschke conjecture}, that {\sl any Blaschke manifold is
isometric to its model CROSS}, had been proved for spherical case by following
scheme [Bes]:

 1) using integral geometry in the space of geodesics one shows that
$vol\,(M^N)\ge vol\,(S^N,can)$ with equality of volumes if and only if $M$ is
isometric to $(S^N,can)$,

 2) on the other hand, one uses topological arguments to show that the {\it
Weinstein integer} $i(M)=vol\,(M^N)/vol\,(S^N,can)$, which has a description
by cogomology of the space of oriented geodesics of given $C_\pi-$manifold, is
actually one.

 The evident analog of Blaschke conjecture for $C_\pi-$manifolds is wrong
already when model CROSS is $S^2$, but unknown for non-spherical case.

 The step 2) is related with {\bf weak Blaschke conjecture} by C.T.Yang, that
{\sl all Blaschke manifolds have right volumes}.
It was proven in [Yan] for complex projective space $CP^n$ and for
$C_\pi-$manifolds homeomorphic to model CROSS in [Rez 1,2].
\vskip 5pt

 In view of facts in 25 the missing link for proving Toponogov's conjecture
is the follows

\proclaim{26 Theorem [RovT]} If a Riemannian manifold $M,\ (\dim M=an)$
has the properties $F_1),\,F_2)$, then its volume is not less then volume of
$KP^n(4)$, $(\dim\Bbb K=a)$ and equality holds if and only if $M$ is
isometric to $KP^n(4)$.\endproclaim

\demo{Proof} We use the scheme [App. D, Bes] with modification by splitting
of Jacobi equation and volume measure along family $\{F(p,\lambda)\}$. The
same proof (for Blaschke manifold with taut geodesics) is given in [Heb].

Let $\mu$ be a volume measure on $M$, the whole measure $V(M)=\int_Md\mu$
of manifold is called {\it volume}. The volumes of projective spaces
$KP^n(4),\ (\Bbb K=\Bbb C,\ \Bbb H,\ \Bbb Ca)$ with standard metrics and
diameter $\frac\pi 2$ and volume of sphere $S^{an-1}$ of radius $1$ will be
denoted by $V(KP^n)$ and $V(S^{an-1})$, their numerical values are known.
The total space of fibration $\pi:UM\to M$ of unit spheres tangent to $M$ is
endowed by canonical metric and measure $\mu_1=\sigma\otimes\mu$, where
$\sigma$ is standard volume measure on $S^{an-1}$.
Thus $V(UM)=V(S^{an-1})V(M)$. The unit {\it geodesical} vector field $Z$ on
$UM$ is defined, -- the projections of the integral curves of $Z$ are
geodesics in $M$, the induced dynamical system $\xi$ on $UM$ is called {\it
geodesical flow}. Note, that the measure $\mu_1$ is invariant under geodesical
flow $\xi$ on manifold $UM$ [Bes].

 Let $\gamma_u,\ (u=(p,\lambda)\in UM)$ be a unit speed geodesic with initial
values $\gamma_u(0)=p,\ \gamma_u'(0)=\lambda$ and $\xi^t$ be a evolution map
of geodesical flow, i.e. $\pi(\xi^t(u))=\gamma_u'(t).$

 Let $f(u,t)$, where $u\in UM,\ t\in R_+$, denotes a volume form on $M$ in
polar coordinates, i.e. at the point $\exp(t\lambda)$ it is true
$d\mu=f(u,t)d\sigma\otimes dt$. This function $f(u,t)$ may be calculated
with the help of Jacobi fields by the following way.
Let $\{\lambda_i\},\ (2\le i\le an)$ be orthonormal basis of subspace
$\lambda^\bot$ (orthogonal complement to $\lambda$) and
$\{Y_i\},\ (2\le i\le an)$ -- Jacobi fields along $\gamma_u$ with initial
values $Y_i(0)=0,\ Y'_i(0)=\lambda_i$. Then for all $t$ it is true:
$f(u,t)=|\det\ Y_2(t)\land\dots\land Y_{an}(t)|$ [Bes]. In our case the
vectors $\lambda_2,\dots,\lambda_a$ are chosen tangent to $d(\lambda)$, and
vectors $\lambda_{a+1},\dots,\lambda_{an}$ -- orthogonal to $d(\lambda)$.
Along geodesic $\gamma_u$ the curvature transformation
$R(*,\gamma'_u),\gamma'_u$ has two invariant subspaces: the tangent space to
the totally geodesic constant curvature submanifold $F(p,\lambda)$ containing
$\gamma_u$, on which it is a multiplication by $4$, and the space orthogonal
to this.

Since the Jacobi fields $Y_2(t),\dots,Y_a(t)$ are tangent to totally
geodesic submanifold $F(p,\lambda)$ with constant curvature 4, they are given
by known formula:
$$Y_i(t)=(\frac 12\sin 2t)\bar\lambda_i,\ (2\le i\le a),$$
where $\bar\lambda_i\ni\lambda_i$ is parallel vector field along $\gamma_u$.
In view of $F_2)$ the vector fields $Y_{a+1}(t),\dots,Y_{an}(t)$ are non-zero
for $0<t<\pi$, they span the normal bundle to $F(p,\lambda)$ along geodesic
$\gamma_u$. Thus,
$$f(u,t)=f_1(u,t)|\frac 12\sin 2t|^{a-1},$$
where function $f_1(u,t)=|\det\ Y_{a+1}(t)\land\dots\land Y_{an}(t)|$ is
positive for $0<t<\pi$. We shall denote
$$\varphi(u,t)=f_1(u,t)^{\frac 1{an-a}},$$
i.e.
$$f(u,t)=\varphi(u,t)^{an-a}|\frac 12\sin 2t|^{a-1}.$$
In particular, for model CROSS $KP^n(4)$ it is true $\ \varphi(u,t)=\sin t$.

\proclaim{27 Lemma}
\ $\pi V^2(M)=\int_{UM}(\int_0^\pi(\int_0^{\pi-x}f(\xi^x(u),t)dt)dx)d\mu_1.$
\endproclaim

\demo {Proof of Lemma 27} Since diameter and injectivity radius of $M$ are
$\frac\pi 2$, then $\forall p\in M$ and a ball $B(p,\frac\pi 2)$ it is true:
$$V(M)=V(B(p,\frac\pi 2))=\int_{U_pM}(\int_0^{\pi/2}f(u,t)dt)d\sigma.$$
With the help of equality $f(-u,t)=f(u,\pi-t),\ (\forall u\in UM)$
[App. D, Bes]), we obtain:
$$\int_{U_pM}(\int_0^{\pi/2}f(-u,t)dt)d\sigma=
\int_{U_pM}(\int_0^{\pi/2}f(u,\pi-t)dt)d\sigma=
\int_{U_pM}(\int_{\pi/2}^\pi f(u,t)dt)d\sigma.\eqno(28)$$
Thus the integration on interval $[0,\pi]$ gives us the double volume
$$2V(M)=\int_{U_pM}(\int_0^\pi f(u,t)dt)d\sigma.\eqno(29)$$
We integrate (29) over $M$
$$2V^2(M)=\int_M(\int_{U_pM}(\int_0^\pi f(u,t)dt)d\sigma)dp=
\int_{UM}(\int_0^\pi f(u,t)dt)d\mu_1$$
and then integrate the last equality over interval $[0,\pi]$ with using the
invariance of measure $\mu_1$ under geodesic flow
$$2\pi V^2(M)=\int_{UM}(\int_0^\pi(\int_0^\pi f(\xi^x(u),t)dt)dx)d\mu_1.
\eqno(30)$$
Since $\xi^x(-u)=-\xi^{\pi-x}(u),\ (u\in UM,\ 0\le x\le\pi)$ and the map
$u\to-u$ is diffeomorphism of manifold $UM$ which preserve measure $\mu_1$,
then in view of $f(-u,t)=f(u,\pi-t)$ it is true
$$\int_{UM}(\int_0^\pi(\int_0^{\pi-x}
\hskip-3pt f(\xi^{\pi-x}(u),\pi-t)dt)dx)d\mu_1=
\int_{UM}(\int_0^\pi(\int^\pi_{\pi-x}\hskip-2pt f(\xi^x(u),t)dt)dx)d\mu_1=$$
$$\int_{UM}(\int_0^\pi (\int_0^{\pi-x}f(\xi^\pi(u),\pi-t)dt)dx)d\mu_1=
\int_{UM}(\int_0^\pi(\int_{\pi-x}^\pi f(\xi^x(u),t)dt)dx)d\mu_1.$$
Thus the integral (30) may be broken onto two equal parts.\qed\enddemo

\proclaim{31 Lemma} For any $u\in UM$ it is true
$$J(u)=
\int_0^\pi(\int_0^{\pi-x}f(\xi^x(u),t)dt)dx\ge\pi\frac{V(KP^n)}{V(S^{an-1})}$$
with equality for only case of
$\ \ K(\gamma_u'(x),y)\equiv 1,\ \ (y\,\bot\,d(\gamma_u'(x)),\ 0\le x\le\pi).$
\endproclaim

\demo{Proof of Lemma 31} We shall use the Holder inequality of order $p=an-a$
$$\int g_1g_2\le(\int g_1^p)^{\frac 1p})
(\int g_2^{\frac p{p-1}})^{\frac{p-1}p},$$
where equality holds if and only if
$g_1^p=\rho g_2^{\frac p{p-1}},\ (\rho=const)$.
In our case the functions will be
$$g_1(t)=\varphi(t)|\frac 12\sin 2t|^\frac{a-1}{an-a},
\ \ g_2(t)=(\sin t)^{an-a-1}|\frac 12\sin 2t|^\frac{(an-a-1)(a-1)}{an-a}$$
and the condition $g_1=\rho g_2^{\frac p{p-1}}$ takes a form:
$\varphi=\rho\sin$. Thus
$$J(u)=\int_0^\pi(\int_x^\pi\varphi^{an-a}(\xi^x(u),y-x)
|\frac 12\sin 2(y-x)|^{a-1}dy)dx\ge$$
$$\frac{(\int_0^\pi\int_x^\pi\varphi(\xi^x(u),y-x)
\sin^{an-a-1}(y-x)\,|\frac 12\sin 2(y-x)|^{a-1}dydx)^{an-a}}
{(\int_0^\pi\int_x^\pi\sin (y-x)^{an-a}\,
|\frac 12\sin 2(y-x)|^{a-1}dydx)^{an-a-1}}.\eqno(32)$$
Since the submanifolds $\{F(p,\lambda)\}$ are totally geodesic, then for
function $\varphi(u,x)$ it is true the inequality (proof is the same, as in
[App. D, Bes])
$$\varphi(\xi^x(u),z)\ge
\varphi(u,x)\varphi(u,x+z)\int_x^{x+z}\frac{dt}{\varphi^2(u,t)}.$$
With the help of Kazdan's inequality [App. E, Bes] with weight function
$$\rho(y-x)=\sin(y-x)^{an-a-1}|\frac 12\sin 2(y-x)|^{a-1}$$
we shall estimate the numerator in (32)
$$\int_0^\pi(\int_x^\pi\varphi(\xi^x(u),y-x)\sin (y-x)^{an-a-1}
\,|\frac 12\sin 2(y-x)|^{a-1}dy)dx\ge$$
$$\int_0^\pi\int_x^\pi\int_x^y
\frac{\varphi(u,x)\varphi(u,y)}{\varphi^2(u,t)}
\sin(y-x)^{an-a-1}|\frac 12\sin 2(y-x)\,|^{a-1}dtdydx\ge$$
$$\int_0^\pi(\int_x^\pi\sin(y-x)^{an-a}\,|\frac 12\sin 2(y-x)|^{a-1}dy)dx=
\beta(a,n).$$
Thus
$$J(u)\ge\beta(a,n)=
\int_0^\pi(\int_x^\pi\sin (y-x)^{an-a}\,|\frac 12\sin 2(y-x)|^{a-1}dy)dx.$$
The equality holds for only case of (see [App. D, Bes])
$$K(\gamma_u'(x),y)=1,\ (0\le x\le\pi,\ y\,\bot\,d(\gamma'_u(x))).$$
We shall show below, that $\beta(a,n)=\pi\frac{V(KP^n)}{V(S^{an-1})}$.
\qed\enddemo

 Continue the proof of Theorem 26. From Lemma 27 and Lemma 31 it follows
$$\pi V^2(M)\ge\int_{UM}\beta(a,n)d\mu_1=V(UM)\beta(a,n).$$
Since $V(UM)=V(M)V(S^{an-1})$, then
$\pi V^2(M)\ge\beta(a,n)V(M)V(S^{an-1})$,
i.e.
$$V(M)\ge\frac 1\pi\beta(a,n)V(S^{an-1}).$$
The equality holds for only case of $M$ being with constant sectional
curvature
$$K(\lambda,y)=1,\ (\lambda\in U_pM,\ y\ \bot\ d(\lambda),\ p\in M).$$
>From above it follows, that our manifold $M$ has positive $\frac 14-$pinched
sectional curvature and by Theorem 18 $M$ is isometric to $KP^n(4)$.
 If we repeat the above for model CROSS $KP^n(4)$, then obtain the numeric
value of $\beta(a,n)$:
$$V(KP^n)=\frac 1\pi\beta (a,n)V(S^{an-1})\ \Rightarrow
\ \beta(a,n)=\pi\frac{V(KP^n)}{V(S^{an-1})}.\qed$$\enddemo

>From the above statement and facts about Blaschke manifolds it follows

\proclaim{33 Theorem [RovT]} Riemannian manifold $M,\ (\dim M=an)$ with
the properties $F_1),\,F_2)$ is isometric to CROSS.\endproclaim

 From Theorem 33 and Theorem 19 it follows the confirmation of
Toponogov's conjecture:

\proclaim{34 Theorem [RovT]} Riemannian manifold $V^m(-\infty,4)$ with
extremal diameter $\frac\pi 2$ is isometric to CROSS.\endproclaim

\proclaim{35 Corollary (diameter rigidity)} A complete, connected, simply
connected Riemannian manifold $M^{2n}$ with sectional curvature $0<K_M\le 4$
and diameter $diam\,(M)=\frac\pi 2$ is isometric to CROSS.\endproclaim

 Theorem 34 and Corollary 35 generalize the Theorem 18 by M.\,Berge.

 Theorem 34 has many corollaries. Below is one of them.
\vskip 5pt

Projective planes $KP^2(4),\ (\Bbb K=\Bbb R,\,\Bbb C,\,\Bbb H,\,\Bbb Ca)$ with
standard metrics are at the same time Riemannian manifolds with dimensions
$2a=2,\,4,\,8,\,16$. Totally geodesic spheres $\{KP^1(4)=S^a(4)\}$ (with
dimension $a$ and constant curvature 4) play the role of {\it straight lines}
in $KP^2(4)$, and {\bf axioms of projective geometry} are true:

{\sl
$P_1$: for all two different points there exists exactly one {\it straight
line}, which connects them,

$P_2$: every  two different {\it straight lines} have intersection at exactly
one point}.

Thus we obtain differential geometrical test of projective planes over
$\Bbb C,\,\Bbb H,\,\Bbb Ca$.

\proclaim{36 Corollary} Assume, that $M^{2a},\ (a>1)$ be a complete
Riemannian manifold with the conditions $P_1,\,P_2$ and {\it straight lines}
are totally geodesic submanifolds isometric to Euclidean sphere $S^a(4)$.
Then $a=2,\,4,\,8$ and $M^{2a}$ is isometric to projective plane $KP^2(4),
\ (\Bbb K=\Bbb C,\,\Bbb H,\,\Bbb Ca)$ with standard metric.\endproclaim

\newpage
\head References\endhead

\enddocument
\end